\documentclass[10pt,a4paper]{article}

\usepackage[T1]{fontenc}
\usepackage[bitstream-charter]{mathdesign}

\usepackage[latin1]{inputenc}							
\usepackage{amsmath}	
\usepackage{latexsym,xypic}								
\usepackage{geometry}
\usepackage{longtable}
\usepackage{array}
\usepackage{fancyvrb,shortvrb}

\usepackage{xcolor}
\definecolor{bl}{rgb}{0.0,0.2,0.6}

\usepackage{sectsty}
\usepackage[compact]{titlesec}
\allsectionsfont{\color{bl}\scshape\selectfont}

\makeatletter							
\def\printtitle{
    {\color{bl} \centering \huge \sc \textbf{\@title}\par}}		
\makeatother							

\title{Computing 2-Dimensional Algebras: Crossed Modules and Cat$^{1}$%
-Algebras  \vspace*{30pt}}

\makeatletter							
\def\printauthor{
    {\color{bl} \centering \large \sc \textbf{\@author}\par}}				
\makeatother							

\author{	Z. Arvasi and A. Odaba\c{s}	\vspace*{10pt}}

\usepackage{fancyhdr}
\usepackage{lastpage}	
\usepackage[runin]{abstract}			
\abslabeldelim{\quad}						%
\newtheorem{theorem}{Theorem}
\newtheorem{definition}[theorem]{Definition}
\newtheorem{example}[theorem]{Example}

\newtheorem{proposition}[theorem]{Proposition}

\begin{document}
\printtitle

\vspace{0.7cm}

\printauthor

\vspace{1cm}

\textbf{Address:} Department of Mathematics and Computer Sciences, Osmangazi University, Art and
Science Faculty, Eskisehir, Turkey.

\vspace{0.5cm}

\textbf{e-mail addresses:} \texttt{zarvasi@ogu.edu.tr, aodabas@ogu.edu.tr}

\vspace{1cm}

\textbf{Abstract:} In this paper, we described the GAP implementation of crossed modules of
commutative algebras and cat$^{1}$-algebras and their equivalence. We include a table of cat$^{1}$-structures on algebras

\vspace{0.5cm}


\vspace{0.5cm}


\section{Introduction}

In 1950 S. MacLane and \ J. H. C. Whitehead, \cite{whitehead}, suggested
that crossed modules modeled homotopy 2-types. Later crossed modules had
been considered as \textquotedblleft 2-dimensional groups\textquotedblright
, \cite{brown1,brown2}. The commutative algebra version of this construction
has been adapted by T. Porter, \cite{arvasi2,porter1}. This algebraic
version is called \textquotedblleft combinatorial algebra
theory\textquotedblright\ which contains potentially important new ideas
(see \cite{arvasi2,arvasi3,arvasi4}.)

A share package XMod, \cite{alp1,alp2}, for the GAP3 computational group
theory language was described by M. Alp and C. Wensley. The 2-dimensional
part of this programme contains functions for computing crossed modules and
cat$^{1}$-groups and their morphisms. This package was rewritten for GAP4,
\cite{gap}, containing functions for crossed square with utility functions,
\cite{alp3}.

In this paper we describe new functions for GAP about computing crossed
module of algebras, cat$^{1}$-algebras and their morphisms by analogy with
\textquotedblleft computational group theory\textquotedblright . The tools
needed are the group algebras in which the group algebra functor $\emph{k}%
(.):\mathbf{Gr}\rightarrow \mathbf{Alg}$ is left adjoint to the unit group
functor $\emph{u}(.):\mathbf{Alg\rightarrow Gr}$.

One of the main result is the GAP implementation of the equivalent
categories \textbf{XModAlg} (crossed modules of algebras) and \textbf{Cat}$%
^{1}$\textbf{-Alg (}cat$^{1}$-algebras\textbf{)} has been shown in this
work. Algorithms of the GAP\ implementations of these structures are deeply
analyzed in the second author's Ph.D. thesis, \cite{aodabas1}.

Unfortunately the GAP\ is not quiet suitable for commutative algebras since
no standard GAP function yet exist for computing semidirect products,
endomorphisms and multiplier algebras. But we manage the endomorphisms and
multiplier algebras. In Section 4, we tabulate, for groups $G$ of order at
most 12 and fields $F$ of order 4, the order of group algebras $FG$; the
number of endomorphisms and idempotent endomorphisms and the number of cat$%
^{1}$-structures on $FG$.

\section{Crossed Modules and Cat$^{1}$- Algebras}

In this section we will present the notion of crossed modules of commutative
algebra and its adaptation to computer environment.

Let \textbf{k} be a fixed commutative ring with $1\neq 0$. From now on, all
\textbf{k}-algebras will be associative and commutative.

A crossed module is a \textbf{k}-algebra morphism $\mathcal{X}:=(\partial
:S\rightarrow R)$ with an action of $R$ on $S$ satisfying%
\[
\begin{array}{ccccccc}
\mathbf{XModAlg1\ :} & \partial (r\cdot s)=r\partial (s), &  &  &  & \mathbf{%
XModAlg2\ :} & \partial (s)\cdot s^{\prime }=ss^{\prime },%
\end{array}%
\]

\noindent for all $s,s^{\prime }\in S$,$\ r\in R$, where $\partial $ is
called the boundary map of $\mathcal{X}$.

A \textbf{k}-algebra homomorphism $\mathcal{X}:=(\partial :S\rightarrow R)$
which satisfy the condition \textbf{XModAlg1} is called a precrossed module.

\bigskip

\textit{Examples :}

\begin{enumerate}
\item[i.] Let $A$ be an \textbf{k}-algebra and $I$ is an ideal of $A$. Then $%
\mathcal{X}=(inc:I\rightarrow A)$ is a crossed module with the
multiplication action of $A$ on $I$. Conversely, we induce an ideal from a
given crossed module. Indeed, for a given crossed module $\mathcal{X}%
=(\partial :S\rightarrow R)$, ${\partial (S)}$ is an ideal of $R$.

\item[ii.] Let $M$ be a $R$-module then $\mathcal{X}=(0:M\rightarrow R)$ is
a crossed module. Conversely, for a crossed module $\mathcal{X}=(\partial
:M\rightarrow R)$, one can get that $Ker\partial $ is a $R/\partial M$%
-module.

\item[iii.] Let $\partial :S\rightarrow R$ be a surjective \textbf{k}%
-algebra homomorphism. Define the action of $R$ on $S$ by $r\cdot s=%
\widetilde{r}s$ where $\widetilde{r}\in \partial ^{-1}(r)$. Then $\mathcal{X}%
=(\partial :S\rightarrow R)$ is a crossed module with the defined action.

\item[iv.] Let $S$ be a \textbf{k}-algebra such that $Ann(S)=0$ or $S^{2}=S$
then $\partial :S\rightarrow M(S)$ is a crossed module, where $M(S)$ is the
algebra of multipliers of $S$ and $\partial $ is the canonical homomorphism.
( see \cite{arvasi4} for this )
\end{enumerate}

A crossed module $\mathcal{X}^{\prime }=(\partial ^{\prime }:S^{\prime
}\rightarrow R^{\prime })$ is a subcrossed module of the crossed module $%
\mathcal{X}=(\partial :S\rightarrow R)$ if $S^{\prime }\leq S$, $R^{\prime
}\leq R$, $\partial ^{\prime }=\partial |_{S^{\prime }}$ and the action of $%
S^{\prime }$ on $R^{\prime }$ is induced by the action of $R$ on $S$.

\begin{definition}
Let $\mathcal{X}=(\partial :S\rightarrow R)$, $\mathcal{X}^{\prime
}=(\partial ^{\prime }:S^{\prime }\rightarrow R^{\prime })$ be (pre) crossed
modules and $\theta :S\rightarrow S^{\prime }$, $\varphi :R\rightarrow
R^{\prime }$ be \textbf{k}-algebra homomorphisms. If%
\[
\begin{array}{ccc}
\varphi \partial =\partial ^{\prime }\theta , &  & \theta (r\cdot s)=\varphi
(r)\cdot \theta (s),%
\end{array}%
\]%
for all $r\in R$ , $s\in S,$ then the pair $(\theta ,\varphi )$ is called a
morphism between $\mathcal{X}=(\partial :S\rightarrow R)$ and $\mathcal{X}%
^{\prime }=(\partial ^{\prime }:S^{\prime }\rightarrow R^{\prime })$.

The conditions can be thought as the commutativity of the following diagrams;

$$ \xymatrix@R=40pt@C=40pt{
  S \ar[d]_{\partial} \ar[r]^{\theta}
                & S^{\prime } \ar[d]^{\partial^{\prime }}  \\
  R  \ar[r]_{\varphi}
                & R^{\prime }             }\ \ \ \
\xymatrix@R=40pt@C=40pt{
  R \times S \ar[d] \ar[r]^{ \varphi \times \theta }
                &  R^{\prime } \times S^{\prime } \ar[d]  \\
  S  \ar[r]_{ \theta } & S^{\prime }       }$$
\end{definition}

Up to this definition we have the category \textbf{XModAlg} of crossed
modules of \textbf{k}-algebras.

Some algebraic structures which are equivalent to crossed modules of
algebras are as follows;

\newpage

\begin{itemize}
\item cat$^{1}$-algebras. (Ellis, \cite{ellis1})

\item simplicial algebras with Moore complex of length 1 (Z. Arvasi and T.
Porter, \cite{arvasi2})

\item algebra-algebroids (Gaffar Musa's Ph.D. thesis, \cite{musa})
\end{itemize}

An alternative algebraic model of $2$-types was defined by Loday in \cite%
{loday} which is called cat$^{1}$-groups. Then algebraic version was given
by Ellis in \cite{ellis1} namely cat$^{1}$-algebras.

Let $A$ and $R$ be \textbf{k}-algebras,%
$$
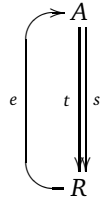
\diagram A \ddto_{t}\ddto<.5ex>^{s} &   \\
\\
      R\uutol^{e}  &     \enddiagram$$
$s,t:A\rightarrow R$ be \textbf{k}-algebra homomorphism and $e:R\rightarrow
A $ is inclusion. If the conditions,

\[
\begin{array}{ccccccc}
\mathbf{Cat1Alg1:} & se=id_{R}=te &  &  &  & \mathbf{Cat1Alg2\ :} & (\ker
s)(\ker t)=\{0_{A}\}%
\end{array}%
\]%
satisfied, then the algebraic system $\mathcal{C}:=(e;t,s:A\rightarrow R)$
is called a cat$^{1}$-algebra. The system which satisfy the condition
Cat1Alg1 is called a precat$^{1}$-algebra. The homomorphisms $s,t$ and $e$
are called as tail, head and embedding homomorphisms, respectively.

\begin{definition}
Let $\mathcal{C}=(e;t,s:A\rightarrow R)$, $\mathcal{C}^{\prime }=(e^{\prime
};t^{\prime },s^{\prime }:A^{\prime }\rightarrow R^{\prime })$ be cat$^{1}$%
-algebras and $\phi :A\rightarrow A^{\prime }$ and $\varphi :R\rightarrow
R^{\prime }$ be \textbf{k}-algebra homomorphisms. If the diagram

$$
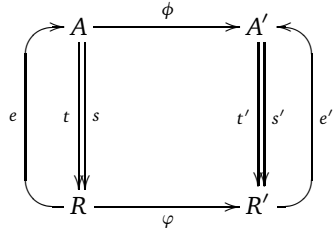
\diagram
A \ddto_{t} \ddto<.5ex>^{s}  %
\rrto^{\phi} &&  A'\ddto_{t'} \ddto<.5ex>^{s'}\\
\\
R \rrto_{\varphi} \uutol^{{e}} && R' \uutor_{e'}
\enddiagram
$$

commutes, then the pair $(\phi ,\varphi )$ is called
a cat$^{1}$-algebra morphism. Consequently we have the category $\mathbf{Cat}^{\mathbf{1}}%
\mathbf{Alg}$ of cat$^{1}$-algebras.
\end{definition}

As indicated in \cite{ellis1},the categories $\mathbf{Cat}^{\mathbf{1}}%
\mathbf{Alg}$ and $\mathbf{XModAlg}$\textbf{\ }are naturally equivalent. For
a given crossed module $\partial :A\rightarrow R$, we have the semidirect
product $R\ltimes A$ thanks to the action of $R$ on $A$. If we define $%
t,s:R\ltimes A\rightarrow R$ and $e:R\rightarrow R\ltimes A$ by%
\[
\begin{tabular}{lllll}
$s(r,a)=r$ &  & $t(r,a)=r+\partial (a)$ &  & $e(r)=(r,0),$%
\end{tabular}%
\]%
respectively, then $\mathcal{C}=(e;t,s:R\ltimes A\rightarrow R)$ is a cat$%
^{1}$-algebra.

Conversely, for a given cat$^{1}$-algebra $\mathcal{C}=(e;t,s:A\rightarrow
R) $, $\partial :ker$ $s\rightarrow $ $R$ is a crossed module, where the
action is conjugate action and $\partial $ is the restriction of $t$ to $%
kers $.

\section{Computer Implementation}

GAP is an open-source system for discrete computational algebra. The system
consists a library of mathematical algorithm implementations, a database
about some algebraic properties of small order groups, vector spaces,
modules, algebras, graphs, codes, designs etc. and some character tables of
these algebraic structures. The system has world wide usage in the area of
education and scientific researches. GAP is free software and user
contributions to the system are supported. These contributions are organized
in a form of GAP packages and are distributed together with the system.
Contributors can submit additional packages for inclusion after a reviewing
process.

The Small Groups library provides access to descriptions of the groups of
small order by Besche, Eick and O'Brien in \cite{bettina}. The groups are
listed up to isomorphism. The library contains the groups those of order at
most 2000 except 1024. There is not similar library for algebras in GAP. To
to this for algebras we will use group algebras as follows:

\vspace{0.3cm}

\noindent \textbf{Group Algebras : }Let \textbf{k }be a field and $G$ be
multiplicative group, finite or infinite. It is well known that the group
algebra \textbf{k}$G$ is an associative \textbf{k}-algebra with the elements
of $G$ as a basis and with multiplication defined distributively using the
group multiplication in $G$. (see \cite{passman} for more knowledge)

Let $\sigma :G\rightarrow H$ be a homomorphism of groups. Define,
\[
\begin{tabular}{cccc}
$\mathbf{k}\sigma :$ & $\mathbf{k}G$ & $\rightarrow $ & $\mathbf{k}H$ \\
& $e_{g}$ & $\mapsto $ & $e_{\sigma (g)}$%
\end{tabular}%
\]
$\mathbf{k}\sigma $ is a homomorphism of group algebras. It has the
properties that $\mathbf{k}id_{G}=id_{\mathbf{k}G}$ and, if $\sigma ^{\prime
}:H\rightarrow J$ is a another group homomorphism $\mathbf{k}(\sigma \sigma
^{\prime })=\mathbf{k}\sigma \mathbf{k}\sigma ^{\prime }$. These facts are
summarized in the following proposition.

\begin{proposition}
$\emph{k}%
(.):\mathbf{Gr}\rightarrow \mathbf{Alg}$ is a functor.
\end{proposition}\label{prop3}

The group algebra functor provides a canonical construction for \textbf{k}%
-algebra from any given group. Conversely, there are at least two canonical
ways of extracting a group from a given \textbf{k}-algebra. One is to forget
the multiplications and take the additive (abelian) group of the algebra;
this gives the forgetful functor $\mathbf{Alg\rightarrow Ab}$.
Alternatively, the subset of the algebra consisting of elements which are
invertible under multiplication forms a subgroup (with the operation of
multiplication, of course) called the group of units of the algebra; this
gives a functor $\emph{u}(.):\mathbf{Alg\rightarrow Gr}$. In general, the
group of units of a non-commutative algebra need not be abelian.

\begin{proposition}
The group algebra functor $\emph{k}%
(.):\mathbf{Gr}\rightarrow \mathbf{Alg}$ is left adjoint to
the unit group functor $\emph{u}(.):\mathbf{Alg\rightarrow Gr}$.
\end{proposition}

\noindent \textbf{Remark : }LAGUNA share package \cite{laguna} extends the
GAP functionality for computations in group algebras. Besides computing some
general properties and attributes of group algebras and their elements. All
of our implementations, we will use group algebras. Thus, we will benefit
from the power of the GAP in the group theory.

\subsection{2-Dimensional Algebras}

In standard GAP Library, there have no functions about crossed modules of
\textbf{k}-algebras and cat$^{1}$-algebras implemented yet. We add
XModAlgByBoundaryAndAction() to construct crossed modules from the given
boundary maps and actions. Additionally we add the functions IsPreXModAlg()
and IsXModAlg() to check the (pre) crossed module structure of the given
systems.

We add the functions XModAlgByIdeal(), XModAlgByModule() and
XModAlgByCentralExtension() to obtain crossed modules from ideals, $R$%
-modules and surjective \textbf{k}-algebra homomorphisms, respectively. Also
we add the function XModAlgByMultipleAlgebra to obtain crossed modules from
given multiplier algebras. Then we add the function XModAlg() which gives
rise to do a choice from the implemented functions, up to input parameters.

Finally, to investigate the properties of the constructed crossed module $CM
$, we add the following functions;

\begin{itemize}
\item The functions Source(CM) and Range(CM), are the source and the range
of the boundary map respectively,

\item The function Boundary(CM), boundary map of the crossed module $CM$,

\item The function XModAlgAction(CM) is the action used in the crossed
module $CM$,

\item The function Size(CM)$,$ is used for calculating the order of the
crossed module $CM$,

\item The function Name(CM)$,$ is used for giving a name for the crossed
module $CM$ by associating the names of source and range algebras,

\item The function Display(CM)$,$ is used for the details of $CM$.
\end{itemize}

The details of these implementations can be found in \cite{aodabas1}.

\begin{example}
In the following example, we construct a crossed module by using the group algebra
$GF_{5}D_{4}$ and its augmentation ideal. Also we show usage of the previous
attributes.

\begin{Verbatim}[frame=single, fontsize=\small, commandchars=\\\{\}]
\textcolor{blue}{gap> A:=GroupRing(GF(5),DihedralGroup(4));}
<algebra-with-one over GF(5), with 2 generators>
\textcolor{blue}{gap> Size(A);}
625
\textcolor{blue}{gap> SetName(A,"GF5[D4]");}
\textcolor{blue}{gap> I:=AugmentationIdeal(A);}
<two-sided ideal in GF5[D4], (2 generators)>
\textcolor{blue}{gap> Size(I);}
125
\textcolor{blue}{gap> SetName(I,"Aug");}
\textcolor{blue}{gap> CM:=XModAlgByIdeal(A,I);}
[Aug->GF5[D(4)]]
\textcolor{blue}{gap> Size(CM);}
[ 125, 625 ]
\textcolor{blue}{gap> Boundary(CM);}
MappingByFunction( Aug, GF5[D4], function( i ) ... end )
\textcolor{blue}{gap> Source(CM);}
Aug
\textcolor{blue}{gap> Range(CM);}
GF5[D4]
\textcolor{blue}{gap> Name(CM);}
"[Aug->GF5[D4]]"
\textcolor{blue}{gap> Display(CM);}
Crossed module [Aug->GF5[D4]] :-
: Source group Aug has generators:
  [ (Z(5)^2)*<identity> of ...+(Z(5)^0)*f1, (Z(5)^2)*<identity> of
  ...+(Z(5)^0)*f2 ]
: Range group GF5[D(4)] has generators:
  [ (Z(5)^0)*<identity> of ..., (Z(5)^0)*f1, (Z(5)^0)*f2 ]
: Boundary homomorphism maps source generators to:
  [ (Z(5)^2)*<identity> of ...+(Z(5)^0)*f1, (Z(5)^2)*<identity> of
  ...+(Z(5)^0)*f2 ]
\end{Verbatim}

\end{example} \label{orn1}

To construct and check the subcrossed modules, we add the functions
SubXModAlg() and IsSubXModAlg().

\begin{example}
In the following GAP session, we constructed a subcrossed module of the crossed module $CM$ constructed in Example 5.

\begin{Verbatim}[frame=single, fontsize=\small, commandchars=\\\{\}]
\textcolor{blue}{gap> eI:=Elements(I);}
\textcolor{blue}{gap> J:=Ideal(I,[eI[4]]);}
<two-sided ideal in Aug, (1 generators)>
\textcolor{blue}{gap> J=I;}
false
\textcolor{blue}{gap> Size(J);}
5
\textcolor{blue}{gap> IsIdeal(I,J);}
true
\textcolor{blue}{gap> IsIdeal(A,J);}
true
\textcolor{blue}{gap> PM:=XModAlg(A,J);}
[Algebra( GF(5),
[ (Z(5)^0)*<identity> of ...+(Z(5)^0)*f1+(Z(5)^2)*f2+(Z(5)^2)*f1*f2
 ] )->GF5[D(4)]]
\textcolor{blue}{gap> Display(PM);}
Crossed module [..->GF5[D4]] :-
: Source group has generators:
  [ (Z(5)^0)*<identity> of ...+(Z(5)^0)*f1+(Z(5)^2)*f2+
  (Z(5)^2)*f1*f2 ]
: Range group GF5[D(4)] has generators:
  [ (Z(5)^0)*<identity> of ..., (Z(5)^0)*f1,
  (Z(5)^0)*f2 ]
: Boundary homomorphism maps source generators to:
  [ (Z(5)^0)*<identity> of ...+(Z(5)^0)*f1+(Z(5)^2)*f2+
  (Z(5)^2)*f1*f2 ]
\textcolor{blue}{gap> IsSubXModAlg(CM,PM);}
true
\end{Verbatim}
\end{example}

We add some new functions about cat$^{1}$-algebras which are not exist in
GAP library. To get cat$^{1}$-algebra structures from given tail, head and
embedding homomorphisms we add the function Cat1AlgByTailHeadEmbedding().
Also, we add the function Cat1AlgByEndomorphisms() to get a structure which
would be a cat$^{1}$-algebra from given two endomorphisms. Then we add the
function Cat1Alg() which makes choice between these two functions up to
input parameters. To check the (pre)cat$^{1}$-algebra structure of this
obtained data, we add the functions IsPreCat1Alg() and IsCat1Alg().

We have the following functions which gives some properties of the
constructed cat$^{1}$-algebra $C$.

\begin{itemize}
\item The functions Source(C) and Range(C), are the source and the range of
the boundary map respectively,

\item The function Boundary(C), boundary map of the cat$^{1}$-algebra $C$,

\item The functions Terminal(C), Initial(C) and Embedding(C) are used for
tail, head and embedding morphisms,

\item The function Size(C)$,$ is used for calculating the order of the $C$,

\item The function Name(C)$,$ is used for giving a name for the cat$^{1}$%
-algebra $C$ by associating the names of source and range algebras.

\item The function Display(C)$,$ is used for the details of $C$.
\end{itemize}

\begin{example}
In the following GAP session, we constructed a cat$^{1}$-algebra on the group algebra obtained from the 4 dimensional Galois field and Klein 4 group.

\begin{Verbatim}[frame=single, fontsize=\small, commandchars=\\\{\}]
\textcolor{blue}{gap> H:=GF(4);}
GF(2^2)
\textcolor{blue}{gap> k4:=Group((1,2),(3,4));}
Group([ (1,2), (3,4) ])
\textcolor{blue}{gap> R:=GroupRing(H,k4);}
<algebra-with-one over GF(2^2), with 2 generators>
\textcolor{blue}{gap> Size(R);}
256
\textcolor{blue}{gap> gR:=GeneratorsOfAlgebra(R);}
[ (Z(2)^0)*(), (Z(2)^0)*(1,2), (Z(2)^0)*(3,4) ]
\textcolor{blue}{gap> f:=AlgebraHomomorphismByImages(R,R,gR,gR);}
[ (Z(2)^0)*(), (Z(2)^0)*(1,2), (Z(2)^0)*(3,4) ] ->
[ (Z(2)^0)*(), (Z(2)^0)*(1,2), (Z(2)^0)*(3,4) ]
\textcolor{blue}{gap> IsAlgebraHomomorphism(f);}
true
\textcolor{blue}{gap> C:=Cat1Alg(f,f,f);}
[AlgebraWithOne( GF(2^2), [ (Z(2)^0)*(1,2), (Z(2)^0)*(3,4)
 ] ) -> AlgebraWithOne( GF(2^2), [ (Z(2)^0)*(1,2),
 (Z(2)^0)*(3,4) ] )]
\textcolor{blue}{gap> IsCat1Alg(C);}
true
\textcolor{blue}{gap> Size(C);}
[ 256, 256 ]
\textcolor{blue}{gap> Display(C);}
Cat1-algebra [..=>..] :-
: source algebra has generators:
  [ (Z(2)^0)*(), (Z(2)^0)*(1,2), (Z(2)^0)*(3,4) ]
:  range algebra has generators:
  [ (Z(2)^0)*(), (Z(2)^0)*(1,2), (Z(2)^0)*(3,4) ]
: tail homomorphism maps source generators to:
  [ (Z(2)^0)*(), (Z(2)^0)*(1,2), (Z(2)^0)*(3,4) ]
: head homomorphism maps source generators to:
  [ (Z(2)^0)*(), (Z(2)^0)*(1,2), (Z(2)^0)*(3,4) ]
: range embedding maps range generators to:
  [ (Z(2)^0)*(), (Z(2)^0)*(1,2), (Z(2)^0)*(3,4) ]
: the kernel is trivial.
\end{Verbatim}

\end{example}

We add the function AllHomsOfAlgebras() to get all homomorphisms between two
algebras, which gives rise to get the functions for constructing cat$^{1}$%
-algebra easier.

The function SubCat1Alg() and IsSubCat1Alg() are added for constructing and
checking the subcat$^{1}$-algebras.

In the following GAP sessions we construct a subcat$^{1}$-algebras by using
these functions.

\begin{Verbatim}[frame=single, fontsize=\small, commandchars=\\\{\}]
\textcolor{blue}{gap> A := GroupRing(GF(2),Group((1,2,3)(4,5)));}
<algebra-with-one over GF(2), with 1 generators>
\textcolor{blue}{gap> R := GroupRing(GF(2),Group((1,2,3)));}
<algebra-with-one over GF(2), with 1 generators>
\textcolor{blue}{gap> f := AllHomsOfAlgebras(A,R);}
[ [ (Z(2)^0)*(1,3,2)(4,5) ] -> [ <zero> of ... ],
  [ (Z(2)^0)*(1,3,2)(4,5) ] -> [ (Z(2)^0)*() ],
  [ (Z(2)^0)*(1,3,2)(4,5) ] -> [ (Z(2)^0)*()+(Z(2)^0)*(1,2,3) ],
  [ (Z(2)^0)*(1,3,2)(4,5) ] ->
    [ (Z(2)^0)*()+(Z(2)^0)*(1,2,3)+(Z(2)^0)*(1,3,2) ],
  [ (Z(2)^0)*(1,3,2)(4,5) ] -> [ (Z(2)^0)*()+(Z(2)^0)*(1,3,2) ],
  [ (Z(2)^0)*(1,3,2)(4,5) ] -> [ (Z(2)^0)*(1,2,3) ],
  [ (Z(2)^0)*(1,3,2)(4,5) ] -> [ (Z(2)^0)*(1,2,3)+(Z(2)^0)*(1,3,2) ],
  [ (Z(2)^0)*(1,3,2)(4,5) ] -> [ (Z(2)^0)*(1,3,2) ] ]
\textcolor{blue}{gap> g := AllHomsOfAlgebras(R,A);}
[ [ (Z(2)^0)*(1,2,3) ] -> [ <zero> of ... ],
  [ (Z(2)^0)*(1,2,3) ] -> [ (Z(2)^0)*() ],
  [ (Z(2)^0)*(1,2,3) ] -> [ (Z(2)^0)*()+(Z(2)^0)*(1,2,3) ],
  [ (Z(2)^0)*(1,2,3) ] -> [ (Z(2)^0)*()+(Z(2)^0)*(1,2,3)+(Z(2)^0)*(1,3,2) ],
  [ (Z(2)^0)*(1,2,3) ] -> [ (Z(2)^0)*()+(Z(2)^0)*(1,3,2) ],
  [ (Z(2)^0)*(1,2,3) ] -> [ (Z(2)^0)*(1,2,3) ],
  [ (Z(2)^0)*(1,2,3) ] -> [ (Z(2)^0)*(1,2,3)+(Z(2)^0)*(1,3,2) ],
  [ (Z(2)^0)*(1,2,3) ] -> [ (Z(2)^0)*(1,3,2) ] ]
\textcolor{blue}{gap> Cat1AlgByTailHeadEmbedding(f[6],f[6],g[8]);}
[AlgebraWithOne( GF(2), [ (Z(2)^0)*(1,2,3)(4,5) ] ) -> AlgebraWithOne( GF(2),
[ (Z(2)^0)*(1,2,3) ] )]
\textcolor{blue}{gap> IsCat1Alg(C);}
true
\textcolor{blue}{gap> Size(C);}
[ 64, 8 ]
\textcolor{blue}{gap> eA := Elements(A);;}
\textcolor{blue}{gap> AA := Subalgebra(A,[eA[1],eA[2],eA[3]]);}
<algebra over GF(2), with 3 generators>
\textcolor{blue}{gap> A = AA;}
false
\textcolor{blue}{gap> eR := Elements(R);;}
\textcolor{blue}{gap> RR := Subalgebra(R,[eR[1],eR[2]]);}
<algebra over GF(2), with 2 generators>
\textcolor{blue}{gap> R = RR;}
false
\textcolor{blue}{gap> CC := SubCat1Alg(C,AA,RR);}
[Algebra( GF(2), [ <zero> of ..., (Z(2)^0)*(), (Z(2)^0)*()+(Z(2)^0)*(4,5)
 ] ) -> Algebra( GF(2), [ <zero> of ..., (Z(2)^0)*() ] )]
\textcolor{blue}{gap> IsCat1Alg(CC);}
true
\textcolor{blue}{gap> Size(CC);}
[ 4, 2 ]
\textcolor{blue}{gap> IsSubCat1Alg(C,CC);}
true
\end{Verbatim}

\subsection{Morphisms of 2-Dimensional Algebras}

First of all, we add the function Make2AlgMorphism() for to define morphisms
of $2$-dimensional algebras such as cat$^{1}$-algebras and crossed modules.
Later we add the function IdentityMapping() to construct the identity
morphisms of given cat$^{1}$-algebra and crossed module.

Let $\mathcal{X}=(\partial :S\rightarrow R)$, $\mathcal{X}^{\prime
}=(\partial ^{\prime }:S^{\prime }\rightarrow R^{\prime })$ be (pre)crossed
modules and $\theta :S\rightarrow S^{\prime }$, $\varphi :R\rightarrow
R^{\prime }$ be algebra homomorphisms. We add the functions
PreXModAlgMorphismByHoms() and XModAlgMorphismByHoms() to construct
precrossed module and crossed module morphisms, respectively, from these two
algebra homomorphisms. One can check the (pre)crossed module morphism
structure of the constructed pair by the added functions
IsPreXModAlgMorphism() and IsXModAlgMorphism().

To get again the pair $(\theta ,\varphi )$ from the constructed crossed
module morphism we add the functions SourceHom() and RangeHom().

Additionally we add the functions Source(), Range() and Kernel() to obtain
the source, range and kernel. On the other hand, we add the functions
IsInjective(), IsSurjecive() and IsBijective() to checking subjectivity and
injectivity of the constructed (pre)crossed module homomorphism.

\begin{example}
In the following GAP session; we will construct a crossed module by using the group algebra obtained from Galois field of order 2 and cyclic group of order 6 and its augmentation ideal. Then we get a crossed module homomorphism by using the implemented functions.

\begin{Verbatim}[frame=single, fontsize=\small, commandchars=\\\{\}]
\textcolor{blue}{gap> A:=GroupRing(GF(2),CyclicGroup(6));}
<algebra-with-one over GF(2), with 2 generators>
\textcolor{blue}{gap> B:=AugmentationIdeal(A);}
<two-sided ideal in <algebra-with-one over GF(2), with 2 generators>
,(dimension 5)>
\textcolor{blue}{gap> CM:=XModAlg(A,B);}
[Algebra( GF(2), [ (Z(2)^0)*<identity> of ...+(Z(2)^0)*f1,
  (Z(2)^0)*f1+(Z(2)^0)*f2, (Z(2)^0)*f2+(Z(2)^0)*f1*f2,
  (Z(2)^0)*f1*f2+(Z(2)^0)*f2^2, (Z(2)^0)*f2^2+(Z(2)^0)*f1*f2^2
 ] )->AlgebraWithOne( GF(2), [ (Z(2)^0)*f1, (Z(2)^0)*f2 ] )]
\textcolor{blue}{gap> SetName(CM,"GF_2C_6");}
\textcolor{blue}{gap> f:=IdentityMapping(CM);}
[[..] => [..]]
\textcolor{blue}{gap> IsPreXModAlgMorphism(f);}
true
\textcolor{blue}{gap> IsXModAlgMorphism(f);}
true
\textcolor{blue}{gap> PM:=Kernel(f);}
[Algebra( GF(2), [], <zero> of ... )->Algebra( GF(2), [],
<zero> of ... )]
\textcolor{blue}{gap> IsXModAlg(PM);}
true
\textcolor{blue}{gap> IsSubXModAlg(CM,PM);}
true
\textcolor{blue}{gap> theta:=SourceHom(f);}
IdentityMapping( <two-sided ideal in <algebra-with-one of dimension
6 over GF(2)>, (dimension 5)> )
\textcolor{blue}{gap> phi:=RangeHom(f);}
IdentityMapping( <algebra-with-one of dimension 6 over GF(2)> )
\textcolor{blue}{gap> IsInjective(f);}
true
\textcolor{blue}{gap> IsSurjective(f);}
true
\textcolor{blue}{gap> IsBijective(f);}
true
\end{Verbatim}

\end{example}

Let $\mathcal{C}=(e;t,s:A\rightarrow R)$, $\mathcal{C}^{\prime }=(e^{\prime
};t^{\prime },s^{\prime }:A^{\prime }\rightarrow R^{\prime })$ be cat$^{1}$%
-algebras and $\phi :A\rightarrow A^{\prime }$ and $\varphi :R\rightarrow
R^{\prime }$ be algebra homomorphisms. To construct a (pre)cat$^{1}$-algebra
morphism we add the functions PreCat1AlgMorphismByHoms() and
Cat1AlgMorphismByHoms(), respectively. One can check the (pre)cat$^{1}$%
-algebra structure of the constructed data, by adding the functions
IsPreCat1AlgMorphism() and IsCat1AlgMorphismMorphism().

We improve the functions about the crossed module morphisms, to that we get
properties of improve given cat$^{1}$-algebra morphisms.

\begin{example}

In the following GAP session, we give an example which shows the construction of a cat$^{1}$-algebra morphism.

\begin{Verbatim}[frame=single, fontsize=\small, commandchars=\\\{\}]
\textcolor{blue}{gap>  A := GroupRing(GF(2),Group(()));}
<algebra-with-one over GF(2), with 1 generators>
\textcolor{blue}{gap> gA := GeneratorsOfAlgebra(A);}
[ (Z(2)^0)*(), (Z(2)^0)*() ]
\textcolor{blue}{gap> m := AlgebraHomomorphismByImages(A,A,gA,gA);}
[ (Z(2)^0)*(), (Z(2)^0)*() ] -> [ (Z(2)^0)*(), (Z(2)^0)*() ]
\textcolor{blue}{gap> C1 := Cat1Alg(m,m,m);}
[AlgebraWithOne( GF(2), [ (Z(2)^0)*() ] ) -> AlgebraWithOne( GF(2),
[ (Z(2)^0)*() ] )]
\textcolor{blue}{gap> B := GroupRing(GF(2),Group((1,2)));}
<algebra-with-one over GF(2), with 1 generators>
\textcolor{blue}{gap> f := AllHomsOfAlgebras(B,A);}
[ [ (Z(2)^0)*(1,2) ] -> [ <zero> of ... ],
  [ (Z(2)^0)*(1,2) ] -> [ (Z(2)^0)*() ] ]
\textcolor{blue}{gap> g := AllHomsOfAlgebras(A,B);}
[ [ (Z(2)^0)*(), (Z(2)^0)*() ] -> [ <zero> of ..., <zero> of ... ],
  [ (Z(2)^0)*(), (Z(2)^0)*() ] -> [ (Z(2)^0)*(), (Z(2)^0)*() ] ]
\textcolor{blue}{gap> C2 := Cat1Alg(f[2],f[2],g[2]);}
[AlgebraWithOne( GF(2), [ (Z(2)^0)*(1,2) ] ) -> AlgebraWithOne( GF(2),
[ (Z(2)^0)*() ] )]
\textcolor{blue}{gap> C1=C2;}
false
\textcolor{blue}{gap> R1:=Source(C1);;}
\textcolor{blue}{gap> R2:=Source(C2);;}
\textcolor{blue}{gap> S1:=Range(C1);;}
\textcolor{blue}{gap> S2:=Range(C2);;}
\textcolor{blue}{gap> gR1:=GeneratorsOfAlgebra(R1);}
[ (Z(2)^0)*(), (Z(2)^0)*() ]
\textcolor{blue}{gap> gR2:=GeneratorsOfAlgebra(R2);}
[ (Z(2)^0)*(), (Z(2)^0)*(1,2) ]
\textcolor{blue}{gap> gS1:=GeneratorsOfAlgebra(S1);}
[ (Z(2)^0)*(), (Z(2)^0)*() ]
\textcolor{blue}{gap> gS2:=GeneratorsOfAlgebra(S2);}
[ (Z(2)^0)*(), (Z(2)^0)*() ]
\textcolor{blue}{gap> im1:=[gR2[1],gR2[1]];}
[ (Z(2)^0)*(), (Z(2)^0)*() ]
\textcolor{blue}{gap> f1:=AlgebraHomomorphismByImages(R1,R2,gR1,im1);}
[ (Z(2)^0)*(), (Z(2)^0)*() ] -> [ (Z(2)^0)*(), (Z(2)^0)*() ]
\textcolor{blue}{gap> im2:=[gS2[1],gS2[1]];}
[ (Z(2)^0)*(), (Z(2)^0)*() ]
\textcolor{blue}{gap> f2:=AlgebraHomomorphismByImages(S1,S2,gS1,im2);}
[ (Z(2)^0)*(), (Z(2)^0)*() ] -> [ (Z(2)^0)*(), (Z(2)^0)*() ]
\textcolor{blue}{gap> m:=Cat1AlgMorphism(C1,C2,f1,f2);}
[[GF(2)_triv=>GF(2)_triv] => [GF(2)_c2=>GF(2)_triv]]
\textcolor{blue}{gap> Display(m);}
Morphism of cat1-algebras :-
: Source = [GF(2)_triv=>GF(2)_triv] with generating sets:
  [ (Z(2)^0)*(), (Z(2)^0)*() ]
  [ (Z(2)^0)*(), (Z(2)^0)*() ]
:  Range = [GF(2)_c2=>GF(2)_triv] with generating sets:
  [ (Z(2)^0)*(), (Z(2)^0)*(1,2) ]
  [ (Z(2)^0)*(), (Z(2)^0)*() ]
: Source Homomorphism maps source generators to:
  [ (Z(2)^0)*(), (Z(2)^0)*() ]
: Range Homomorphism maps range generators to:
  [ (Z(2)^0)*(), (Z(2)^0)*() ]
\textcolor{blue}{gap> Image2dAlgMapping(m);}
[GF(3)_c2^3=>GF(3)_c2^3]
\textcolor{blue}{gap> IsSurjective(m);}
false
\textcolor{blue}{gap> IsInjective(m);}
true
\textcolor{blue}{gap> IsBijective(m);}
false
\end{Verbatim}

\end{example}

\subsection{Naturally Equivalence}

By using the natural equivalence of categories of crossed modules and cat$%
^{1}$-algebras, we add the functions XModAlgByCat1Alg() and
Cat1AlgByXModAlg() which constructs crossed modules and cat$^{1}$-algebras
from the given cat$^{1}$-algebras and crossed modules, respectively.

\begin{example}

In the following GAP session, we get a cat$^{1}$-algebra from a crossed module.

\begin{Verbatim}[frame=single, fontsize=\small, commandchars=\\\{\}]
\textcolor{blue}{gap> R:=GroupRing(GF(3),CyclicGroup(2));}
<algebra-with-one over GF(3), with 1 generators>
\textcolor{blue}{gap> I:=AugmentationIdeal(R);}
<two-sided ideal in <algebra-with-one over GF(3), with 1 generators>,
  (1 generators)>
\textcolor{blue}{gap> CM:=XModAlgByIdeal(R,I);}
[Algebra( GF(3), [ (Z(3))*<identity> of ...+(Z(3)^0)*f1
 ] )->AlgebraWithOne( GF(3), [ (Z(3)^0)*f1 ] )]
\textcolor{blue}{gap> IsXModAlg(CM);}
true
\textcolor{blue}{gap> C:=Cat1AlgByXModAlg(CM);}
[AlgebraWithOne( GF(3), [ (Z(3)^0)*f1 ] ) IX Algebra( GF(3),
[ (Z(3))*<identity> of ...+(Z(3)^0)*f1 ] ) -> AlgebraWithOne( GF(3),
[ (Z(3)^0)*f1 ] )]
\textcolor{blue}{gap> IsCat1Alg(C);}
true
\textcolor{blue}{gap> SM:=XModAlgByCat1Alg(C);}
[..->..]
\textcolor{blue}{gap> SM=CM;}
true
\end{Verbatim}
\end{example}

Since all these operations linked to the functions \texttt{Cat1Alg} and
\texttt{XModAlg}, then all of them can be done by using these two functions.

\begin{example} In the following GAP session, we get a crossed module from a  cat$^{1}$-algebra.

\begin{Verbatim}[frame=single, fontsize=\small, commandchars=\\\{\}]
\textcolor{blue}{gap> StructureDescription(Group((1,2),(2,3)));}
"S3"
\textcolor{blue}{gap> A := GroupRing(GF(4),Group((1,2),(2,3)));}
<algebra-with-one over GF(2^2), with 2 generators>
\textcolor{blue}{gap> Size(A);}
4096
\textcolor{blue}{gap> m := AlgebraHomomorphismByImages(A,A,gA,gA);}
[ (Z(2)^0)*(), (Z(2)^0)*(1,2), (Z(2)^0)*(2,3) ] ->
[ (Z(2)^0)*(), (Z(2)^0)*(1,2), (Z(2)^0)*(2,3) ]
\textcolor{blue}{gap> C := Cat1AlgByEndomorphisms(m,m);}
[AlgebraWithOne( GF(2^2), [ (Z(2)^0)*(1,2), (Z(2)^0)*(2,3)
 ] ) -> AlgebraWithOne( GF(2^2), [ (Z(2)^0)*(1,2), (Z(2)^0)*(2,3) ] )]
\textcolor{blue}{gap> CM:=XModAlg(C);}
[Algebra( GF(2^2), [], <zero> of ... )->AlgebraWithOne( GF(2^2),
[ (Z(2)^0)*(1,2), (Z(2)^0)*(2,3) ] )]
\textcolor{blue}{gap> IsXModAlg(CM);}
true
\textcolor{blue}{gap> CC:=Cat1Alg(CM);}
[AlgebraWithOne( GF(2^2), [ (Z(2)^0)*(1,2), (Z(2)^0)*(2,3)
 ] ) -> AlgebraWithOne( GF(2^2), [ (Z(2)^0)*(1,2), (Z(2)^0)*(2,3) ] )]
\textcolor{blue}{gap> CC=C;}
true
\end{Verbatim}

\end{example}

\section{Table of Cat$^{1}$-Structures}

Every idempotent endomorphism pair of a algebra is a candidate for
constructing a cat$^{1}$-algebra structure. At this vein, we found all
endomorphisms of the algebras which were given in the following table. Then
we obtain the idempotent endomorphisms and construct all cat$^{1}$-algebras
defined from the given group algebras by using the functions given in
section 3.

For each group $G$ and field $F$ we list size of group algebra $FG$; the
size of the set endomorphisms $End(FG)$ on $FG$; the size of the set $IE(FG)$
of idempotents in $End(FG)$; and the size of the set of all cat$^{1}$%
-structures $C(FG)$ on $FG$.

\begin{longtable}{ccccccc}
\hline
Field & Gap id & Group & $\left\vert FG\right\vert $ & $\left\vert \text{%
End(F}G\text{)}\right\vert $ & $\left\vert \text{IE(F}G\text{)}\right\vert $
& $\left\vert C(FG)\right\vert $ \\ \hline\hline
GF(2) & [1,1] & I & 2 & 2 & 2 & 1 \\
GF(3) & [1,1] & I & 3 & 2 & 2 & 1 \\
GF(4) & [1,1] & I & 4 & 2 & 2 & 1 \\
GF(2) & [2,1] & C2 & 4 & 3 & 3 & 2 \\
GF(3) & [2,1] & C2 & 9 & 9 & 6 & 3 \\
GF(4) & [2,1] & C2 & 16 & 5 & 3 & 2 \\
GF(2) & [3,1] & C3 & 8 & 8 & 5 & 1 \\
GF(3) & [3,1] & C3 & 27 & 10 & 3 & 1 \\
GF(4) & [3,1] & C3 & 64 & 64 & 23 & 7 \\
GF(2) & [4,1] & C4 & 16 & 9 & 3 & 1 \\
GF(3) & [4,1] & C4 & 81 & 45 & 18 & 3 \\
GF(4) & [4,1] & C4 & 256 & 65 & 3 & 1 \\
GF(2) & [4,2] & C2$\times $C2 & 16 & 65 & 15 & 13 \\
GF(3) & [4,2] & C2$\times $C2 & 81 & 625 & 104 & 25 \\
GF(4) & [4,2] & C2$\times $C2 & 256 & 4097 & 83 & 81 \\
GF(2) & [5,1] & C5 & 32 & 12 & 5 & 1 \\
GF(3) & [5,1] & C5 & 243 & 12 & 5 & 1 \\
GF(4) & [5,1] & C5 & 1024 & 72 & 21 & 5 \\
GF(2) & [6,1] & S3 & 64 & 51 & 23 & 2 \\
GF(3) & [6,1] & S3 & 729 & 201 & 37 & 7 \\
GF(4) & [6,1] & S3 & 4096 & 485 & 63 & 2 \\
GF(2) & [6,2] & C6 & 64 & 39 & 14 & 4 \\
GF(3) & [6,2] & C6 & 729 & 361 & 29 & 13 \\
GF(4) & [6,2] & C6 & 4096 & 2197 & 168 & 44 \\
GF(2) & [7,1] & C7 & 128 & 128 & 25 & 7 \\
GF(3) & [7,1] & C7 & 2187 & 16 & 5 & 1 \\
GF(4) & [7,1] & C7 & 16584 & 128 & 25 & 7 \\
GF(2) & [8,1] & C8 & 256 & 129 & 3 & 1 \\
GF(3) & [8,1] & C8 & 6561 & 6561 & 468 & 39 \\
GF(4) & [8,1] & C8 & 65536 & 16385 & 3 & 1 \\
GF(2) & [8,2] & C4$\times $c2 & 256 & 8193 & 131 & 65 \\
GF(2) & [8,3] & D8 & 256 & 2305 & 67 & 1 \\
GF(2) & [8,4] & Q8 & 256 & 1793 & 3 & 1 \\
GF(2) & [8,5] & C2$\times $C2$\times $C2   & 256 & 2657423 & 723  & 87 \\
GF(2) & [9,1] & C9 & 512 & 80 & 17 & 1 \\
GF(3) & [9,1] & C9 & 19683 & 6562 & 3 & 1 \\
GF(4) & [9,1] & C9 & 262144 & 6400 & 419 & 49 \\
GF(2) & [9,2] & C3$\times $c3 & 512 & 20000 & 809 & 73 \\
GF(2) & [10,1] & D10 & 1024 & 471 & 63 & 2 \\
GF(2) & [10,2] & C10 & 1024 & 243 & 26 & 4 \\
GF(2) & [11,1] & C11 & 2048 & 24 & 5 & 1 \\
GF(2) & [12,1] & C3:C4 & 4096 & 1881 & 167 & 9 \\
GF(2) & [12,2] & C12 & 4096 & 1737 & 74 & 1 \\
GF(2) & [12,3] & A4 & 4096 & 2210 & 179 & 1 \\
GF(2) & [12,4] & D12 & 4096 & 38545 & 1583 & 117 %
\end{longtable}

\section*{Acknowledgements}

This works was partially supported by T\"{U}B\.{I}TAK (The Scientific and
Technical Research Council of Turkey).

Project Number : 107T542


\begin{thebibliography}{99}
\bibitem{alp1} \textsc{Alp M. and C.D. Wensley} Crossed Modules and Cat%
${{}^1}$%
-groups in GAP, \emph{version 1 Manual for the XMOD share package for GAP3},
(1996).

\bibitem{alp3} \textsc{Alp M. and C.D. Wensley} Crossed Modules and Cat%
${{}^1}$%
-groups in GAP, \emph{version 2.23 Manual for the XMOD share package for GAP4%
}, (2013).

\bibitem{alp2} \textsc{Alp M. and C.D. Wensley} Enumeration of cat$^{1}$%
-groups of low order, \emph{Int. J. Algebra Comput.}, 10, 407-424, (2000).

\bibitem{arvasi2} \textsc{Arvasi Z. and Porter T.} Simplicial and Crossed
Resolutions of Commutative Algebras, \emph{Journal of Algebra}, 181,
426-448, (1996).

\bibitem{arvasi3} \textsc{Arvasi Z. and Porter T.} Freeness Conditions for
2-Crossed Modules of Commutative Algebras, \emph{Applied Categorical
Structures}, 6, 455-471, (1998).

\bibitem{arvasi4} \textsc{Arvasi Z. and Ege U.} Annihilators, multipliers
and crossed modules, \emph{Applied Categorical Structures}, 11, 487--506
(2003).

\bibitem{bettina} \textsc{Besche H.U., Eick B. and O'Brien E.A.} A
millennium project: constructing Small Groups, \emph{\ Internat. J. Algebra
Comput.}, 12, 623 - 644 (2002).

\bibitem{laguna} \textsc{Bovdi V. and et al} Lie AlGebras and UNits of group
Algebras in GAP, \emph{version 3.64 Manual for the LAGUNA share package for
GAP4}, (2013).

\bibitem{brown1} \textsc{Brown R.} Higher Dimensional Group Teory, Low
Dimensional Topology, \emph{London Math. Soc. Lecture Note Series}, 48,
215-238, (1982).

\bibitem{brown2} \textsc{Brown R.} From groups to groupoids: a brief survey,%
\emph{Bull. London Math. Soc.}, 19, 113-134, (1987).

\bibitem{ellis1} \textsc{Ellis G.J.} Higher Dimensional Crossed Modules of
Algebras. \ \emph{J.Pure Appl. Algebra 52}, 277-282, (1988).

\bibitem{gap} \textsc{\ GAP - Groups, Algortihms, and Programming, Version 4}
Lehrstuhl D f\"{u}r Mathematik, RWTH Aachen Germany and School of
Mathematical and Computational Sciences, \emph{U. St. Andrews, Scotland},
(1997).

\bibitem{loday} \textsc{Loday J.L.} Spaces with finitely many non-trivial
homotopy groups., \emph{J. App. Algebra}, 24 , 179--202, (1982).

\bibitem{musa} \textsc{Mosa G. H.} Higher dimensional algebroids and crossed
complexes. \ \emph{Ph.D. Thesis}, \ University of Wales, Bangor (U.K.)
1991-94 \ (1987).

\bibitem{aodabas1} \textsc{Odabas A.} Crossed Modules of Algebras with GAP.
\ \emph{Ph.D. Thesis}, \ Osmangazi University. \ (2009).

\bibitem{passman} \textsc{Passman D.S.} The Algebraic Structure of Group
Rings. \emph{A Wiley-Interscience Publication}, \ (1977).

\bibitem{porter1} \textsc{Porter T.} Some categorical results in the theory
of crossed modules in commutative algebras. \emph{J. Algebra}, 109, 415-429,
(1987).

\bibitem{whitehead} \textsc{Whitehead J.H.C.} Combinatorial Homotopy II,
\emph{Bull. American Math. Society}, 55, 453-456, (1949).
\end{thebibliography}
\end{document}